\documentclass[10pt]{article}
\usepackage{amsmath} 
\usepackage{amssymb}
\usepackage{latexsym}
\usepackage{xcolor}
\usepackage{fancyhdr}
\allowdisplaybreaks 

 \usepackage{comment}

\def\XXint#1#2#3{{\setbox0=\hbox{$#1{#2#3}{\int}$} 
\vcenter{\hbox{$#2#3$}}\kern-.5\wd0}}   

 \numberwithin{equation}{section}
\newtheorem{theorem}[equation]{Theorem}
\newtheorem{proposition}[equation]{Proposition}
\newtheorem{definition}[equation]{Definition}
\newtheorem{remark}[equation]{Remark}
\newtheorem{lemma}[equation]{Lemma}

\title{
A limiting case  of a theorem of C. Miranda for   layer potentials in Schauder spaces} 
 
\author{  
Massimo Lanza de Cristoforis
\\
Dipartimento di Matematica `Tullio Levi-Civita', 
\\
Universit\`a degli Studi di Padova, 
\\
Via Trieste 63, Padova 35121, 
Italy. 
\\
E-mail: mldc@math.unipd.it   }

\date{\ }

\begin{document}
 
 \maketitle

\noindent
{\bf Abstract:}  The aim of this paper is to prove a theorem of C.~Miranda for the single and double layer potential corresponding to the fundamental solution of a  second order differential operator with constant coefficients in Schauder spaces in the limiting case in which the open set is of class $C^{m,1}$ and  the densities are of class $C^{m-1,1}$ for the single layer potential and of class $C^{m,1}$ for the double layer potential
for some nonzero natural number $m$. The treatment of the limiting case requires generalized Schauder spaces.

 \vspace{\baselineskip}



\section{Introduction} In this paper, we consider the single and  double layer potential associated to the fundamental solution of a second order differential operator with constant coefficients. Unless otherwise specified, we assume throughout the paper that
\[
n\in {\mathbb{N}}\setminus\{0,1\}\,,
\]
where ${\mathbb{N}}$ denotes the set of natural numbers including $0$. Let $\alpha\in]0,1]$, $m\in {\mathbb{N}} $. Let $\Omega$ be a bounded open subset of ${\mathbb{R}}^{n}$ of class $C^{m,\alpha}$. Let
\[
\Omega^-\equiv {\mathbb{R}}^{n}\setminus\overline{\Omega}
\]
denote the exterior of $\Omega$.  For the notation and standard properties of the (generalized) Schauder spaces and sets of class $C^{m,\alpha}$ we refer for example  to Dondi and the author  \cite[\S 2]{DoLa17} and to the reference \cite[\S 2.11, 2.13]{DaLaMu21} of Dalla Riva, the author and Musolino. Let $\nu \equiv (\nu_{l})_{l=1,\dots,n}$ denote the external unit normal to $\partial\Omega$.  Let $N_{2}$ denote the number of multi-indexes $\gamma\in {\mathbb{N}}^{n}$ with $|\gamma|\leq 2$. For each 
\begin{equation}
\label{mldc_introd0}
{\mathbf{a}}\equiv (a_{\gamma})_{|\gamma|\leq 2}\in {\mathbb{C}}^{N_{2}}\,, 
\end{equation}
we set 
\[
a^{(2)}\equiv (a_{lj} )_{l,j=1,\dots,n}\qquad
a^{(1)}\equiv (a_{j})_{j=1,\dots,n}\qquad
a\equiv a_{0} \,,
\]
with $a_{lj} \equiv 2^{-1}a_{e_{l}+e_{j}}$ for $j\neq l$, $a_{jj} \equiv
 a_{e_{j}+e_{j}}$,
and $a_{j}\equiv a_{e_{j}}$, where $\{e_{j}:\,j=1,\dots,n\}$  is the canonical basis of ${\mathbb{R}}^{n}$. We note that the matrix $a^{(2)}$ is symmetric. 
Then we assume that 
  ${\mathbf{a}}\in  {\mathbb{C}}^{N_{2}}$ satisfies the following ellipticity assumption
\begin{equation}
\label{mldc_ellip}
\inf_{
\xi\in {\mathbb{R}}^{n}, |\xi|=1
}{\mathrm{Re}}\,\left\{
 \sum_{|\gamma|=2}a_{\gamma}\xi^{\gamma}\right\} >0\,,
\end{equation}
and we consider  the case in which
\begin{equation}
\label{mldc_symr}
a_{lj} \in {\mathbb{R}}\qquad\forall  l,j=1,\dots,n\,.
\end{equation}
Then we introduce the operators
\begin{eqnarray*}
P[{\mathbf{a}},D]u&\equiv&\sum_{l,j=1}^{n}\partial_{x_{l}}(a_{lj}\partial_{x_{j}}u)
+
\sum_{l=1}^{n}a_{l}\partial_{x_{l}}u+au\,,
\\
B_{\Omega}^{*}v&\equiv&\sum_{l,j=1}^{n} \overline{a}_{jl}\nu_{l}\partial_{x_{j}}v
-\sum_{l=1}^{n}\nu_{l}\overline{a}_{l}v\,,
\end{eqnarray*}
for all $u,v\in C^{2}(\overline{\Omega})$,  a fundamental solution $S_{{\mathbf{a}} }$ of $P[{\mathbf{a}},D]$,   the single layer potential
\[
v_\Omega[S_{{\mathbf{a}}}   ,\mu](x) =\int_{\partial\Omega}
S_{ {\mathbf{a}} }(x-y)\mu (y)\,d\sigma_{y} \qquad\forall x\in {\mathbb{R}}^{n}\,,
\]
and the double layer potential 
\begin{eqnarray}
\label{mldc_introd3}
\lefteqn{
w_\Omega[{\mathbf{a}},S_{{\mathbf{a}}}   ,\mu](x) \equiv 
\int_{\partial\Omega}\mu (y)\overline{B^{*}_{\Omega,y}}\left(S_{{\mathbf{a}}}(x-y)\right)
\,d\sigma_{y}
}
\\  \nonumber
&&
\qquad
=-\int_{\partial\Omega}\mu(y)\sum_{l,j=1}^{n} a_{jl}\nu_{l}(y)\frac{\partial S_{ {\mathbf{a}} } }{\partial x_{j}}(x-y)\,d\sigma_{y}
\\  \nonumber
&&
\qquad\quad
-\int_{\partial\Omega}\mu(y)\sum_{l=1}^{n}\nu_{l}(y)a_{l}
S_{ {\mathbf{a}} }(x-y)\,d\sigma_{y} \qquad\forall x\in {\mathbb{R}}^{n}\,,
\end{eqnarray}
where the density or moment $\mu$ is a function  from $\partial\Omega$ to ${\mathbb{C}}$. Here the subscript $y$ of $\overline{B^{*}_{\Omega,y}}$ means that we are taking $y$ as variable of the differential operator $\overline{B^{*}_{\Omega,y}}$. 

It is well known that if $\mu$ is $\alpha$-H\"{o}lder continuous, and $\Omega$ is of class $C^{1,\alpha}$, then the restriction 
$v_\Omega[S_{{\mathbf{a}}}   ,\mu]_{|\Omega}$ has a continuous extension 
$v_\Omega^+[S_{{\mathbf{a}}}   ,\mu]$ to $\overline{\Omega}$, that 
$v_\Omega[S_{{\mathbf{a}}}   ,\mu]_{|\Omega^-}$ has a continuous extension 
$v_\Omega^-[S_{{\mathbf{a}}}   ,\mu]$ to $\overline{\Omega^-}$, that
$w_\Omega[{\mathbf{a}},S_{{\mathbf{a}}}   ,\mu]_{|\Omega}$ has a continuous extension 
$w_\Omega^+[{\mathbf{a}},S_{{\mathbf{a}}}   ,\mu]$ to $\overline{\Omega}$ and  that 
$w_\Omega[{\mathbf{a}},S_{{\mathbf{a}}}   ,\mu]_{|\Omega^-}$ has a continuous extension 
$w_\Omega^-[{\mathbf{a}},S_{{\mathbf{a}}}   ,\mu]$ to $\overline{\Omega^-}$. 

In a remarkable paper, Miranda \cite{Mi65} has considered the case of homogeneous differential operators and
has shown that 
if $\alpha\in]0,1[$,   $m\geq 1$,  $\Omega$ is a bounded open subset of ${\mathbb{R}}^n$ of class $C^{m,\alpha}$, then 
$v_\Omega^+[S_{{\mathbf{a}}}   ,\cdot]$ defines a   linear and continuous operator from
$C^{m-1,\alpha}(\partial\Omega)$ to $C^{m,\alpha}(\overline{\Omega})$ and that 
$w_\Omega^+[{\mathbf{a}},S_{{\mathbf{a}}}   ,\cdot]$ defines a   linear and continuous operator from
$C^{m,\alpha}(\partial\Omega)$ to $C^{m,\alpha}(\overline{\Omega})$ and that a corresponding result holds in the exterior of $\Omega$. 
Miranda  \cite{Mi65} has also considered elliptic operators of order greater or equal to $2$. 

 In the present paper, we plan to consider case the limiting case $\alpha=1$, a case that requires the introduction of the  
 generalized Schauder space $C^{m,\omega_{1}(\cdot)}(\overline{\Omega})$ 
  of the functions of class $C^{m}(\overline{\Omega})$ such that the  $m$-th order  partial derivatives satisfy a generalized $\omega_{1}$-H\"{o}lder condition with
\begin{equation}
\label{mldc_omth-1}
\omega_{1}(r)
\equiv
\left\{
\begin{array}{ll}
0 &r=0\,,
\\
r |\ln r | &r\in]0,r_{1}]\,,
\\
r_{1} |\ln r_{1} | & r\in ]r_{1},+\infty[\,,
\end{array}
\right.
\end{equation}
where
$
r_{1}\equiv e^{-1}
$ and we show that $v_\Omega^+[S_{{\mathbf{a}}}   ,\cdot]$ defines a   linear and continuous operator from
$C^{m-1,1}(\partial\Omega)$ to $C^{m,\omega_{1}(\cdot)}(\overline{\Omega})$ (cf.~Theorems \ref{mldc_thm:sla010om}, \ref{mldc_thm:dslam1mom})
and that 
$w_\Omega^+[{\mathbf{a}},S_{{\mathbf{a}}}   ,\cdot]$ defines a   linear and continuous operator from
$C^{m,1}(\partial\Omega)$ to $C^{m,\omega_{1}(\cdot)}(\overline{\Omega})$    and that a corresponding result holds in the exterior of $\Omega$ (cf.~Theorems \ref{mldc_thm:dlay111om}, \ref{mldc_thm:dslam1mom}). 

The paper is organized as follows. Section \ref{mldc_sec:tecprel}   is a section of preliminaries and notation. In Section \ref{mldc_sec:prelfuso}, we introduce some properties on the fundamental solution $S_{{\mathbf{a}} }$ that we need. In Section \ref{mldc_sec:1miraslay} we   prove our Theorem \ref{mldc_thm:sla010om} on the single layer potential in case $m=1$. In Section \ref{mldc_sec:1miradolay} we   prove our Theorem \ref{mldc_thm:dlay111om} on the double layer potential in case $m=1$.   In section \ref{mldc_sec:miradslay} we   prove our Theorem \ref{mldc_thm:dslam1mom} on the single and double layer potential in case $m\geq 1$. 
 
\section{Preliminaries and notation}\label{mldc_sec:tecprel}
If $X$ and $Y$, $Z$ are normed spaces, then ${\mathcal{L}}(X,Y)$ denotes the space of linear and continuous maps from $X$ to $Y$ and ${\mathcal{L}}^{(2)}(X\times Y, Z)$ 
denotes the space of bilinear and continuous maps from $X\times Y$ to $Z$ with their usual operator norm (cf.~\textit{e.g.}, \cite[pp.~16, 621]{DaLaMu21}). 
$|A|$ denotes the operator norm of a matrix $A$ with real (or complex) entries, 
       $A^{t}$ denotes the transpose matrix of $A$. $\delta_{l,j}$ denotes the Kronecker\index{Kronecker symbol}  symbol. Namely,  $\delta_{l,j}=1$ if $l=j$, $\delta_{l,j}=0$ if $l\neq j$, with $l,j\in {\mathbb{N}}$.  The symbol
$| \cdot|$ denotes the Euclidean modulus   in
${\mathbb{R}}^{n}$ or in ${\mathbb{C}}$. For all $r\in]0,+\infty[$, $ x\in{\mathbb{R}}^{n}$, 
$x_{j}$ denotes the $j$-th coordinate of $x$, and  
 ${\mathbb{B}}_{n}( x,r)$ denotes the ball $\{y\in{\mathbb{R}}^{n}:\, |x- y|<r\}$.  If ${\mathbb{D}}$ is a subset of $ {\mathbb{R}}^n$, 
then we set
\[
B({\mathbb{D}})\equiv\left\{
f\in {\mathbb{C}}^{\mathbb{D}}:\,f\ \text{is\ bounded}
\right\}
\,,\quad
\|f\|_{B({\mathbb{D}})}\equiv\sup_{\mathbb{D}}|f|\qquad\forall f\in B({\mathbb{D}})\,.
\]
Then $C^0({\mathbb{D}})$ denotes the set of continuous functions from ${\mathbb{D}}$ to ${\mathbb{C}}$ and we introduce the subspace
$
C^0_b({\mathbb{D}})\equiv C^0({\mathbb{D}})\cap B({\mathbb{D}})
$
of $B({\mathbb{D}})$.  Let $\omega$ be a function from $[0,+\infty[$ to itself such that
\begin{eqnarray}
\nonumber
&&\qquad\qquad\omega(0)=0,\qquad \omega(r)>0\qquad\forall r\in]0,+\infty[\,,
\\
\label{mldc_om}
&&\qquad\qquad\omega\ {\text{is\   increasing,}}\ \lim_{r\to 0^{+}}\omega(r)=0\,,
\\
\nonumber
&&\qquad\qquad{\text{and}}\ \sup_{(a,t)\in[1,+\infty[\times]0,+\infty[}
\frac{\omega(at)}{a\omega(t)}<+\infty\,.
\end{eqnarray}
Here `$\omega$ is increasing' means that 
$\omega(r_1)\leq \omega(r_2)$ whenever $r_1$, $r_2\in [0,+\infty[$ and $r_1<r_2$.
If $f$ is a function from a subset ${\mathbb{D}}$ of ${\mathbb{R}}^n$   to ${\mathbb{C}}$,  then we denote by   $|f:{\mathbb{D}}|_{\omega (\cdot)}$  the $\omega(\cdot)$-H\"older constant  of $f$, which is delivered by the formula   
\[
|f:{\mathbb{D}}|_{\omega (\cdot)
}
\equiv
\sup\left\{
\frac{|f( x )-f( y)|}{\omega(| x- y|)
}: x, y\in {\mathbb{D}} ,  x\neq
 y\right\}\,.
\]        
If $|f:{\mathbb{D}}|_{\omega(\cdot)}<+\infty$, we say that $f$ is $\omega(\cdot)$-H\"{o}lder continuous. Sometimes, we simply write $|f|_{\omega(\cdot)}$  
instead of $|f:{\mathbb{D}}|_{\omega(\cdot)}$. The
subset of $C^{0}({\mathbb{D}} ) $  whose
functions  are
$\omega(\cdot)$-H\"{o}lder continuous    is denoted  by  $C^{0,\omega(\cdot)} ({\mathbb{D}})$
and $|f:{\mathbb{D}}|_{\omega(\cdot)}$ is a semi-norm on $C^{0,\omega(\cdot)} ({\mathbb{D}})$.  
Then we consider the space  $C^{0,\omega(\cdot)}_{b}({\mathbb{D}} ) \equiv C^{0,\omega(\cdot)} ({\mathbb{D}} )\cap B({\mathbb{D}} ) $ with the norm \[
\|f\|_{ C^{0,\omega(\cdot)}_{b}({\mathbb{D}} ) }\equiv \sup_{x\in {\mathbb{D}} }|f(x)|+|f|_{\omega(\cdot)}\qquad\forall f\in C^{0,\omega(\cdot)}_{b}({\mathbb{D}} )\,.
\] 
\begin{remark}
\label{mldc_rem:om4}
Let $\omega$ be as in (\ref{mldc_om}). 
Let ${\mathbb{D}}$ be a   subset of ${\mathbb{R}}^{n}$. Let $f$ be a bounded function from $ {\mathbb{D}}$ to ${\mathbb{C}}$, $a\in]0,+\infty[$.  Then,
\[
\sup_{x,y\in {\mathbb{D}},\ |x-y|\geq a}\frac{|f(x)-f(y)|}{\omega(|x-y|)}
\leq \frac{2}{\omega(a)} \sup_{{\mathbb{D}}}|f|\,.
\]
\end{remark}
In the case in which $\omega(\cdot)$ is the function 
$r^{\alpha}$ for some fixed $\alpha\in]0,1]$, a so-called H\"{o}lder exponent, we simply write $|\cdot:{\mathbb{D}}|_{\alpha}$ instead of
$|\cdot:{\mathbb{D}}|_{r^{\alpha}}$, $C^{0,\alpha} ({\mathbb{D}})$ instead of $C^{0,r^{\alpha}} ({\mathbb{D}})$, $C^{0,\alpha}_{b}({\mathbb{D}})$ instead of $C^{0,r^{\alpha}}_{b} ({\mathbb{D}})$, and we say that $f$ is $\alpha$-H\"{o}lder continuous provided that 
$|f:{\mathbb{D}}|_{\alpha}<+\infty$.   The function $\omega_{1}(\cdot) $ of (\ref{mldc_omth-1}) is concave and satisfies   condition (\ref{mldc_om}).
 We also note that if ${\mathbb{D}}\subseteq {\mathbb{R}}^n$, then the continuous embedding
\begin{equation}\label{mldc_eq:0t0om0t}
C^{0, 1 }_b({\mathbb{D}})\subseteq 
C^{0,\omega_1(\cdot)}_b({\mathbb{D}})\subseteq 
C^{0,\theta'}_b({\mathbb{D}})
\end{equation}
holds for all $\theta'\in ]0,1[$. For the standard properties of the spaces of H\"{o}lder or Lipschitz continuous functions, we refer    to \cite[\S 2]{DoLa17}, \cite[\S 2.6]{DaLaMu21}.\par

Let $\Omega$ be an open subset of ${\mathbb{R}}^n$. The space of $m$ times continuously 
differentiable complex-valued functions on $\Omega$ is denoted by 
$C^{m}(\Omega,{\mathbb{C}})$, or more simply by $C^{m}(\Omega)$. 
Let $f\in  C^{m}(\Omega) $. Then   $Df$ denotes the Jacobian matrix of $f$. 
Let  $\eta\equiv
(\eta_{1},\dots ,\eta_{n})\in{\mathbb{N}}^{n}$, $|\eta |\equiv
\eta_{1}+\dots +\eta_{n}  $. Then $D^{\eta} f$ denotes
$\frac{\partial^{|\eta|}f}{\partial
x_{1}^{\eta_{1}}\dots\partial x_{n}^{\eta_{n}}}$.    The
subspace of $C^{m}(\Omega )$ of those functions $f$ whose derivatives $D^{\eta }f$ of
order $|\eta |\leq m$ can be extended with continuity to 
$\overline{\Omega}$  is  denoted $C^{m}(
\overline{\Omega})$. 
The
subspace of $C^{m}(\overline{\Omega} ) $  whose
functions have $m$-th order derivatives that are
H\"{o}lder continuous  with exponent  $\alpha\in
]0,1]$ is denoted $C^{m,\alpha} (\overline{\Omega})$ and the
subspace of $C^{m}(\overline{\Omega} ) $  whose
functions have $m$-th order derivatives that are
$\omega(\cdot)$-H\"{o}lder continuous    is denoted $C^{m,\omega(\cdot)} (\overline{\Omega})$. 

The subspace of $C^{m}(\overline{\Omega} ) $ of those functions $f$ such that  $f_{|\overline{(\Omega\cap{\mathbb{B}}_{n}(0,r))}}$ belongs to $
C^{m,\omega}(\overline{(\Omega\cap{\mathbb{B}}_{n}(0,r))})$ for all $r\in]0,+\infty[$ is denoted $C^{m,\omega}_{{\mathrm{loc}}}(\overline{\Omega} ) $. \par
Now let $\Omega $ be a bounded
open subset of  ${\mathbb{R}}^{n}$. Then $C^{m}(\overline{\Omega} )$,  $C^{m,\omega(\cdot) }(\overline{\Omega})$ with $\omega$ as in (\ref{mldc_om}) and 
  $C^{m,\alpha }(\overline{\Omega})$ with $\alpha\in]0,1]$   are endowed with their usual norm and are well known to be 
Banach spaces  (cf.~\textit{e.g.},  \cite[\S 2]{DoLa17}, Dalla Riva, the author and Musolino  \cite[\S 2.11]{DaLaMu21}).

  For the (classical) definition of the H\"{o}lder and Schauder spaces    $C^{m,\alpha}(\partial\Omega)$, $C^{m,\omega(\cdot)}(\partial\Omega)$
on the boundary $\partial\Omega$ of an open set $\Omega$ for some $m\in{\mathbb{N}}$, $\alpha\in ]0,1]$, we refer for example to Dondi and the author \cite[\S 2]{DoLa17}, 
    Dalla Riva, the author and Musolino  \cite[\S  2.20]{DaLaMu21}.

 The space of real valued functions of class $C^\infty$ with compact support in an open set $\Omega$ of ${\mathbb{R}}^n$ is denoted ${\mathcal{D}}(\Omega)$. Then its dual ${\mathcal{D}}'(\Omega)$ is known to be the space of distributions in $\Omega$. The support of a function is denoted by the abbreviation `${\mathrm{supp}}$'.  

Morever, we retain the standard notation for the Lebesgue spaces $L^p$ for $p\in [1,+\infty]$ (cf.~\textit{e.g.}, Folland \cite[Chapt.~6]{Fo95}, \cite[\S 2.1]{DaLaMu21}) and
$m_n$ denotes the
$n$ dimensional Lebesgue measure. If $n\in {\mathbb{N}}\setminus\{0\}$, $m\in {\mathbb{N}}$, $h\in {\mathbb{R}}$, $\alpha\in ]0,1]$, then we set 
\[
{\mathcal{K}}^{m,\alpha}_h \equiv\biggl\{
k\in C^{m,\alpha}_{ {\mathrm{loc}}}({\mathbb{R}}^n\setminus\{0\}):\, k\ {\text{is\ positively\ homogeneous\ of \ degree}}\ h
\biggr\}\,,
\]
where $C^{m,\alpha}_{ {\mathrm{loc}}}({\mathbb{R}}^n\setminus\{0\})$ denotes the set of functions of 
$C^{m}({\mathbb{R}}^n\setminus\{0\})$ whose restriction to $\overline{\Omega}$ is of class $C^{m,\alpha}(\overline{\Omega})$ for all bounded open subsets $\Omega$ of ${\mathbb{R}}^n$ such that
$\overline{\Omega}\subseteq {\mathbb{R}}^n\setminus\{0\}$
and we set
\[
\|k\|_{ {\mathcal{K}}^{m,\alpha}_h}\equiv \|k\|_{C^{m,\alpha}(\partial{\mathbb{B}}_n(0,1))}\qquad\forall k\in {\mathcal{K}}^{m,\alpha}_h\,.
\]
We can easily verify that $ \left({\mathcal{K}}^{m,\alpha}_h , \|\cdot\|_{ {\mathcal{K}}^{m,\alpha}_h}\right)$ is a Banach space and we  consider the closed subspaces
\begin{eqnarray*}
{\mathcal{K}}^{m,\alpha}_{h;o}&\equiv& \left\{
k\in {\mathcal{K}}^{m,\alpha}_h:\,k\ {\text{is\ odd}} 
\right\}\,,
\\ \nonumber
{\mathcal{K}}^{m,\alpha}_{h;e,0}&\equiv& \left\{
k\in {\mathcal{K}}^{m,\alpha}_h:\,k\ {\text{is\ even}}, \int_{\partial{\mathbb{B}}_n(0,1)}k\,d\sigma=0
\right\}
\end{eqnarray*}
of ${\mathcal{K}}^{m,\alpha}_h$.  We now introduce the following extension of a classical result of Miranda \cite{Mi65} (see also \cite[Thm.~4.17 (ii)]{DaLaMu21}), who has considered the case of domains of class $C^{1,\alpha}$  and of densities $\mu \in C^{0,\alpha}(\partial\Omega)$ for $\alpha\in]0,1[$ to the limiting case in which $\alpha=1$. 
 For a proof, we refer to Dalla Riva, the author and Musolino \cite{DaLaMu24a}.
 \begin{theorem}\label{mldc_thm:mirandao01}
Let $\Omega$ be a bounded open subset of ${\mathbb{R}}^n$ of class $C^{1,1}$. Then the following statements hold.
\begin{enumerate}
\item[(i)] For each $(k,\mu)\in {\mathcal{K}}^{1,1}_{-(n-1);o}\times C^{0,1}(\partial\Omega)$, the map
\[
  \int_{\partial\Omega} k(x-y) \mu(y)\,d\sigma_{y}\qquad \forall x\in\Omega 
\]
can be extended to a unique  $\omega_1(\cdot)$-H\"{o}lder continuous function $K[k,\mu]^+$ on $\overline{\Omega}$. Moreover, the map from ${\mathcal{K}}^{1,1}_{-(n-1);o}\times C^{0,1}(\partial\Omega)$ to $C^{0,\omega_1(\cdot)}(\overline{\Omega})$ that  takes $(k,\mu)$ to $K[k,\mu]^+$ is bilinear and continuous.
\item[(ii)] Let $r\in]0,+\infty[$ be such that $\overline{\Omega}\subseteq {\mathbb{B}}_n(0,r)$. Then for each   $(k,\mu)\in {\mathcal{K}}^{1,1}_{-(n-1);o}\times C^{0,1}(\partial\Omega)$  the map
\[
  \int_{\partial\Omega} k(x-y) \mu(y)\,d\sigma_{y}\qquad \forall x\in
 {\mathbb{R}}^n\setminus\overline{\Omega}\,,
\]
can be extended to a unique continuous function $K[k,\mu]^-$ on ${\mathbb{R}}^n\setminus \Omega$ such that 
the restriction $K[k,\mu]^-_{|\overline{{\mathbb{B}}_n(0,r)}\setminus\Omega}$ is $\alpha$-H\"{o}lder continuous. Moreover,  the map from ${\mathcal{K}}^{1,1}_{-(n-1);o}\times C^{0,1}(\partial\Omega)$ to $C^{0,\omega_1(\cdot)}( \overline{{\mathbb{B}}_n(0,r)}\setminus\Omega)$ that takes $(k,\mu)$ to $K[k,\mu]^-_{|\overline{{\mathbb{B}}_n(0,r)}\setminus\Omega}$ is bilinear and continuous.
\end{enumerate}
\end{theorem}

\section{Preliminaries on the fundamental solution}
  \label{mldc_sec:prelfuso}

In order to analyze the volume potential, we need some more information on the fundamental solution $S_{ {\mathbf{a}} } $.  To do so, we introduce the fundamental solution $S_{n}$ of the Laplace operator. Namely, we set
\[
S_{n}(x)\equiv
\left\{
\begin{array}{lll}
\frac{1}{s_{n}}\ln  |x| \qquad &   \forall x\in 
{\mathbb{R}}^{n}\setminus\{0\},\quad & {\mathrm{if}}\ n=2\,,
\\
\frac{1}{(2-n)s_{n}}|x|^{2-n}\qquad &   \forall x\in 
{\mathbb{R}}^{n}\setminus\{0\},\quad & {\mathrm{if}}\ n>2\,,
\end{array}
\right.
\]
where $s_{n}$ denotes the $(n-1)$ dimensional measure of 
$\partial{\mathbb{B}}_{n}(0,1)$ and
we follow a formulation of Dalla Riva \cite[Thm.~5.2, 5.3]{Da13} and Dalla Riva, Morais and Musolino \cite[Thm.~5.5]{DaMoMu13}, that we state  as in paper  \cite[Cor.~4.2]{DoLa17} with Dondi (see also John~\cite{Jo55}, Miranda~\cite{Mi65} for homogeneous operators, and Mitrea and Mitrea~\cite[p.~203]{MitMit13}).   
\begin{proposition}
 \label{mldc_prop:ourfs} 
Let ${\mathbf{a}}$ be as in (\ref{mldc_introd0}), (\ref{mldc_ellip}), (\ref{mldc_symr}). 
Let $S_{ {\mathbf{a}} }$ be a fundamental solution of $P[{\mathbf{a}},D]$. 
Then there exist an invertible matrix $T\in M_{n}({\mathbb{R}})$ such that
\begin{equation}
\label{mldc_prop:ourfs0}
a^{(2)}=TT^{t}\,,
\end{equation}
 a real analytic function $A_{1}$ from $\partial{\mathbb{B}}_{n}(0,1)\times{\mathbb{R}}$ to ${\mathbb{C}}$ such that
 $A_{1}(\cdot,0)$ is odd,    $b_{0}\in {\mathbb{C}}$, a real analytic function $B_{1}$ from ${\mathbb{R}}^{n}$ to ${\mathbb{C}}$ such that $B_{1}(0)=0$, and a real analytic function $C $ from ${\mathbb{R}}^{n}$ to ${\mathbb{C}}$ such that
\begin{eqnarray*}
\lefteqn{S_{ {\mathbf{a}} }(x)
= 
\frac{1}{\sqrt{\det a^{(2)} }}S_{n}(T^{-1}x)
}
\\ \nonumber
&&\qquad\qquad
+|x|^{3-n}A_{1}(\frac{x}{|x|},|x|)
 +(B_{1}(x)+b_{0}(1-\delta_{2,n}))\ln  |x|+C(x)\,,
\end{eqnarray*}
for all $x\in {\mathbb{R}}^{n}\setminus\{0\}$,
 and such that both $b_{0}$ and $B_{1}$   equal zero
if $n$ is odd. Moreover, 
 \[
 \frac{1}{\sqrt{\det a^{(2)} }}S_{n}(T^{-1}x) 
 \]
is a fundamental solution for the principal part
  of $P[{\mathbf{a}},D]$.
\end{proposition}
In particular for the statement that $A_{1}(\cdot,0)$ is odd, we refer to
Dalla Riva, Morais and Musolino \cite[Thm.~5.5, (32)]{DaMoMu13}, where $A_{1}(\cdot,0)$ coincides with ${\mathbf{f}}_1({\mathbf{a}},\cdot)$ in that paper. Here we note that a function $A$ from $(\partial{\mathbb{B}}_{n}(0,1))\times{\mathbb{R}}$ to ${\mathbb{C}}$ is said to be real analytic provided that it has a real analytic extension   to an open neighbourhood of $(\partial{\mathbb{B}}_{n}(0,1))\times{\mathbb{R}}$ in 
${\mathbb{R}}^{n+1}$. Then we have the following elementary lemma    (cf.~\textit{e.g.},  \cite[Lem.~4.2]{La22d}).		 
\begin{lemma}\label{mldc_lem:anexsph}
 Let $n\in {\mathbb{N}}\setminus\{0,1\}$. A function $A$ from   $(\partial{\mathbb{B}}_{n}(0,1))\times{\mathbb{R}}$ to ${\mathbb{C}}$ is  real analytic if and only if the function $\tilde{A}$ from $({\mathbb{R}}^n\setminus\{0\}) \times{\mathbb{R}}$ defined by 
\[
\tilde{A}(x,r)\equiv A(\frac{x}{|x|},r)\qquad\forall (x,r)\in ({\mathbb{R}}^n\setminus\{0\}) \times{\mathbb{R}}
\]
is real analytic.
 \end{lemma}
 
Then  one can prove the following formula for the Jacobian of the fundamental solution (see  Dondi and the author	 	 \cite[Lem.~4.3,  (4.8) and the following 2 lines]{DoLa17}). Here one should remember that $A_1(\cdot,0)$ is odd and that $b_0=0$ if $n$ is odd). 
\begin{proposition}
\label{mldc_prop:grafun}
 Let ${\mathbf{a}}$ be as in (\ref{mldc_introd0}), (\ref{mldc_ellip}), (\ref{mldc_symr}). Let $T\in M_{n}({\mathbb{R}})$  be as in (\ref{mldc_prop:ourfs0}). Let $S_{ {\mathbf{a}} }$ be a fundamental solution of $P[{\mathbf{a}},D]$. Let  $B_{1}$, $C$
 be as in Proposition \ref{mldc_prop:ourfs}. 
  Then there exists a real analytic function $A_{2}\equiv (A_{2,j})_{j=1,\dots,n}$ from $\partial{\mathbb{B}}_{n}(0,1)\times{\mathbb{R}}$ to ${\mathbb{C}}^{n}$ such that
\begin{eqnarray*}
\lefteqn{
DS_{ {\mathbf{a}} }(x)=\frac{1}{ s_{n}\sqrt{\det a^{(2)} } }
|T^{-1}x|^{-n}x^{t}(a^{(2)})^{-1} 
}
\\ \nonumber
&&\qquad\qquad
+|x|^{2-n}A_{2}(\frac{x}{|x|},|x|)+DB_{1}(x)\ln |x|+DC(x)
\end{eqnarray*}
for all $x\in {\mathbb{R}}^{n}\setminus\{0\}$. 
Moreover,   $A_2(\cdot,0)$ is even.
\end{proposition}
Next we introduce some notation. If $X$ and $Y$ are subsets of ${\mathbb{R}}^n$, then the symbol
\[
{\mathbb{D}}_{X\times Y}\equiv  \left\{
(x,y)\in X\times Y:\,x=y
\right\}
\] denotes the diagonal set of $X\times Y$  
 and we introduce the following class of   `potential type' kernels  (see also paper \cite{DoLa17} of the author and Dondi,  where such classes have been introduced in a form that generalizes those of Giraud \cite{Gi34}, Gegelia \cite{Ge67}, 
 Kupradze, Gegelia, Basheleishvili and 
 Burchuladze \cite[Chap.~IV]{KuGeBaBu79}).	 
\begin{definition}
Let $X$, $Y\subseteq {\mathbb{R}}^n$. Let $s_1$, $s_2$, $s_3\in {\mathbb{R}}$. We denote by the symbol ${\mathcal{K}}_{s_1, s_2, s_3} (X\times Y)$ the set of continuous functions $K$ from $(X\times Y)\setminus {\mathbb{D}}_{X\times Y}$ to ${\mathbb{C}}$ such that
 \begin{eqnarray*}
\lefteqn{
\|K\|_{  {\mathcal{K}}_{ s_1, s_2, s_3  }(X\times Y)  }
\equiv
\sup\biggl\{\biggr.
|x-y|^{ s_{1} }\vert K(x,y)\vert :\,(x,y)\in X\times Y, x\neq y
\biggl.\biggr\}
}
\\ \nonumber
&&\qquad\qquad\qquad
+\sup\biggl\{\biggr.
\frac{|x'-y|^{s_{2}}}{|x'-x''|^{s_{3}}}
\vert   K(x',y)- K(x'',y)  \vert :\,
\\ \nonumber
&&\qquad\qquad\qquad 
x',x''\in X, x'\neq x'', y\in Y\setminus {\mathbb{B}}_n(x',2|x'-x''|)
\biggl.\biggr\}<+\infty\,.
\end{eqnarray*}
\end{definition}
 
\section{A limiting case of a theorem of C.~Miranda for the single layer potential in case $m=1$}\label{mldc_sec:1miraslay}
If $h\in C^0({\mathbb{R}}^n\setminus\{0\})$, we find convenient to set
\begin{eqnarray}
\label{mldc_eq:1miraslay}
&&K[h,\mu](x)\equiv\int_{\partial\Omega}h(x-y)\mu(y)\,d\sigma_y\qquad\forall x\in{\mathbb{R}}^n\setminus\partial\Omega\,, 
\\ \nonumber
&&K^+[h,\mu] \equiv K[h,\mu]_{|\Omega}\,,
\qquad
K^-[h,\mu] \equiv K[h,\mu]_{|{\mathbb{R}}^n\setminus\overline{\Omega}}\,,
\end{eqnarray}
for all $\mu\in L^\infty(\partial\Omega)$. In those cases in which $K^+[h,\mu]$ can be extended by continuity to $\overline{\Omega}$, we still denote the extension by the symbol $K^+[h,\mu]$ and in those cases in which $K^-[h,\mu]$ can be extended by continuity to ${\mathbb{R}}^n\setminus \Omega$, we still denote the extension by the symbol $K^-[h,\mu]$. 
Then we   consider the following preliminary technical statement on the single layer potential that corresponds to the fundamental solution of a nonhomogeneous   second order elliptic differential operator with constant coefficients.
 
\begin{proposition}\label{mldc_prop:pdersa010om}
 Let ${\mathbf{a}}$ be as in (\ref{mldc_introd0}), (\ref{mldc_ellip}), (\ref{mldc_symr}).  Let $S_{ {\mathbf{a}} }$ be a fundamental solution of $P[{\mathbf{a}},D]$.  Let $\Omega$ be a bounded open subset of ${\mathbb{R}}^n$ of class $C^{1,1}$. Let $j\in\{1,\dots,n\}$, 
\[
K\left[\frac{\partial S_{ {\mathbf{a}} }}{\partial x_j},\mu\right](x)\equiv
\int_{\partial\Omega}\frac{\partial S_{ {\mathbf{a}} }}{\partial x_j}(x-y)\mu(y)\,d\sigma_y
\qquad\forall x\in {\mathbb{R}}^n\setminus\partial\Omega\,,
\]
for all $\mu\in C^0(\partial\Omega)$. 
Then the following statements hold.
\begin{enumerate}
\item[(i)] For each $\mu\in C^{0,1}(\partial\Omega)$, the map $K\left[\frac{\partial S_{ {\mathbf{a}} }}{\partial x_j},\mu\right]_{|\Omega}$ can be extended to a unique $\omega_1(\cdot)$-H\"{o}lder continuous function $K^+\left[\frac{\partial S_{ {\mathbf{a}} }}{\partial x_j},\mu\right]$ on $\overline{\Omega}$. Moreover, the map from $C^{0,1}(\partial\Omega)$ to $C^{0,\omega_1(\cdot)}(\overline{\Omega})$ that takes $\mu$ to $K^+\left[\frac{\partial S_{ {\mathbf{a}} }}{\partial x_j},\mu\right]$ is linear and continuous.
\item[(ii)]  Let $r\in]0,+\infty[$ be such that $\overline{\Omega}\subseteq {\mathbb{B}}_n(0,r)$. For each $\mu\in C^{0,1}(\partial\Omega)$, the map $K\left[\frac{\partial S_{ {\mathbf{a}} }}{\partial x_j},\mu\right]_{|{\mathbb{R}}^n\setminus\overline{\Omega}}$ can be extended to   continuous function $K^-\left[\frac{\partial S_{ {\mathbf{a}} }}{\partial x_j},\mu\right]$ on ${\mathbb{R}}^n\setminus\Omega$.  Moreover, the map from $C^{0,1}(\partial\Omega)$ to $C^{0,\omega_1(\cdot)}(\overline{{\mathbb{B}}_n(0,r)}\setminus\Omega)$ that takes $\mu$ to $K^-\left[\frac{\partial S_{ {\mathbf{a}} }}{\partial x_j},\mu\right]_{|\overline{{\mathbb{B}}_n(0,r)}\setminus\Omega}$ is linear and continuous.
\end{enumerate}
\end{proposition}
 {\bf Proof.} Let $r\in]0,+\infty[$ be such that $\overline{\Omega}\subseteq {\mathbb{B}}_n(0,r)$.  Let
 \begin{eqnarray*}\nonumber
 &&J_j(x)\equiv \frac{1}{ s_{n}\sqrt{\det a^{(2)} } }
|T^{-1}x|^{-n}\left(x^{t}(a^{(2)})^{-1} \right)_j\,,
\\
&& 
 k_j(x)\equiv \frac{\partial S_{ {\mathbf{a}} }}{\partial x_j}-\frac{1}{ s_{n}\sqrt{\det a^{(2)} } }
|T^{-1}x|^{-n}\left(x^{t}(a^{(2)})^{-1} \right)_j
 \end{eqnarray*}
 for all $x\in {\mathbb{R}}^{n}\setminus\{0\}$. Since $J_j$ is odd and positively homogeneous of degreee $-(n-1)$, Theorem \ref{mldc_thm:mirandao01} implies that
 $K[J_j,\mu]_{|\Omega}$ can be extended by continuity to $\overline{\Omega}$, that
 $K[J_j,\mu]_{|{\mathbb{R}}^n\setminus\overline{\Omega}}$ can be extended by continuity to ${\mathbb{R}}^n\setminus\Omega$ for all $\mu\in C^{0,1}(\partial\Omega)$ and that
 \begin{eqnarray}\label{mldc_prop:pdersa010om2}
&&K^+[J_j,\cdot]\in {\mathcal{L}}\left(C^{0,1}(\partial\Omega),C^{0,\omega_1(\cdot)}(\overline{\Omega})\right)\,,
\\	\nonumber
&&K^-[J_j,\cdot]_{|\overline{{\mathbb{B}}_n(0,r)}\setminus\Omega}\in {\mathcal{L}}\left(C^{0,1}(\partial\Omega),C^{0,\omega_1(\cdot)}(\overline{{\mathbb{B}}_n(0,r)}\setminus\Omega)\right)\,.
\end{eqnarray}
Thus it suffices to show that statements (i) and (ii) are satisfied by the integral with kernel $k_j(x-y)$. Since $\partial\Omega$ is upper $(n-1)$-Ahlfors regular with respect to ${\mathbb{R}}^n$ (cf. \cite[Prop.~6.5]{La24}) and the singularity of $k_j(x-\cdot)$ is weak for all $x\in {\mathbb{R}}^n$, we have
 \begin{eqnarray}\label{mldc_prop:pdersa010om3}
&&K^+[k_j,\cdot]\in {\mathcal{L}}\left(C^{0,1}(\partial\Omega),C^{0}(\overline{\Omega})\right)\,,
\\	\nonumber
&&K^-[k_j,\cdot]_{|\overline{{\mathbb{B}}_n(0,r)}\setminus\Omega}\in {\mathcal{L}}\left(C^{0,1}(\partial\Omega),C^{0}(\overline{{\mathbb{B}}_n(0,r)}\setminus\Omega)\right)\,.
\end{eqnarray}
(cf.~\cite[Prop.~4.1]{La22a}, \cite[Prop.~4.3]{La22b}). By the classical Theorem of differentiation of integrals depending on a parameter, we have
\[
\frac{\partial}{\partial x_l}\int_{\partial\Omega}k_j(x-y)\mu(y)\,d\sigma_y
=
\int_{\partial\Omega}\frac{\partial}{\partial x_l}k_j(x-y)\mu(y)\,d\sigma_y
\qquad\forall x\in {\mathbb{R}}^n\setminus\partial\Omega
\]
for all $l\in\{1,\dots,n\}$.
Thus it suffices to show that 
\begin{eqnarray}\label{mldc_prop:pdersa010om4}
&&K^+\left[\frac{\partial}{\partial x_l}k_j,\cdot\right]\in {\mathcal{L}}\left(C^{0,1}(\partial\Omega),C^{0}(\overline{\Omega})\right)\,,
\\	\nonumber
&&K^-\left[\frac{\partial}{\partial x_l}k_j,\cdot\right]_{|\overline{{\mathbb{B}}_n(0,r)}\setminus\Omega}\in {\mathcal{L}}\left(C^{0,1}(\partial\Omega),C^{0}(\overline{{\mathbb{B}}_n(0,r)}\setminus\Omega)\right)\,, 
\end{eqnarray} 
(see the notation in (\ref{mldc_eq:1miraslay})).
Indeed, the following continuous embeddings 
\begin{eqnarray}\nonumber
&&C^{1}(\overline{\Omega})\subseteq C^{0,1}(\overline{\Omega})
\subseteq
C^{0,\omega_1(\cdot)}(\overline{\Omega})\,,
\\ \nonumber
&&
C^{1}(\overline{ {\mathbb{B}}_n(0,r)}\setminus\Omega)\subseteq C^{0,1}(\overline{{\mathbb{B}}_n(0,r)}\setminus\Omega)
\subseteq
C^{0,\omega_1(\cdot)}(\overline{{\mathbb{B}}_n(0,r)}\setminus\Omega)\,,
\end{eqnarray}
hold true (cf.~(\ref{mldc_eq:0t0om0t}), \cite[Lem.~2.28]{DaLaMu21}). Let  $(\xi_1,\dots,\xi_n,r)$ denote  the
variable of $A_2$. 
Then by classical differentiation rules and by the fundamental theorem of calculus, we have
\begin{eqnarray}\label{mldc_prop:pdersa010om5}
\lefteqn{
\frac{\partial k_j}{\partial x_l}(x)=(2-n)|x|^{1-n}\frac{x_l}{|x|}A_{2,j}(\frac{x}{|x|},|x|)
}
\\ \nonumber
&&\qquad
+|x|^{2-n}\sum_{s=1}^n\frac{\partial A_{2,j}}{\partial\xi_s}(\frac{x}{|x|},|x|)
\left(
\frac{\delta_{sl}}{|x|}-\frac{x_sx_l}{|x|^3}
\right)
+|x|^{2-n}\frac{\partial A_{2,j}}{\partial r}(\frac{x}{|x|},|x|)\frac{x_l}{|x|}
\\ \nonumber
&&\qquad
+\frac{\partial^2 B_1}{\partial x_l\partial x_j}(x)\ln |x|+\frac{\partial B_1}{\partial x_j}(x)\frac{x_l}{|x|^2}+\frac{\partial^2 C}{\partial x_l\partial x_j}(x)
\\ \nonumber
&&\qquad
=J_{jl,1}(x)
+
(2-n)|x|^{2-n}\frac{x_l}{|x|}\int_0^1\frac{\partial A_{2,j}}{\partial r}(\frac{x}{|x|},\tau |x|)\,d\tau
\\ \nonumber
&&\qquad
+J_{jl,2}(x)
+|x|^{3-n}\sum_{s=1}^n\int_0^1\frac{\partial^2 A_{2,j}}{\partial\xi_s\partial r}(\frac{x}{|x|},\tau|x|)\,d\tau
\left(
\frac{\delta_{sl}}{|x|}-\frac{x_sx_l}{|x|^3}
\right)
\\ \nonumber
&&\qquad
+|x|^{2-n}\frac{\partial A_{2,j}}{\partial r}(\frac{x}{|x|}, |x|)\frac{x_l}{|x|}
\\ \nonumber
&&\qquad
+\frac{\partial^2 B_1}{\partial x_l\partial x_j}(x)\ln |x|+\frac{\partial^2 C}{\partial x_l\partial x_j}(x)
\\ \nonumber
&&\qquad
+J_{jl,3}(x)+\int_0^1\sum_{s=1}^n\frac{\partial^2 B_1}{\partial x_j\partial x_s}(\tau x)\,d\tau\frac{x_s x_l}{|x|^2}
\qquad\forall x\in {\mathbb{R}}^{n}\setminus\{0\}\,,
\end{eqnarray}
where
\begin{eqnarray}\label{mldc_prop:pdersa010om6}
&&J_{jl,1}(x)\equiv (2-n)|x|^{1-n}\frac{x_l}{|x|}A_{2,j}(\frac{x}{|x|},0)
\\ \nonumber
&&J_{jl,2}(x)\equiv |x|^{2-n}\sum_{s=1}^n\frac{\partial A_{2,j}}{\partial\xi_s}(\frac{x}{|x|},0)\left(
\frac{\delta_{sl}}{|x|}-\frac{x_sx_l}{|x|^3}
\right)
\\ \nonumber
&&J_{jl,3}(x)\equiv
\frac{\partial B_1}{\partial x_j}(0)\frac{x_l}{|x|^2}
\qquad\forall x\in {\mathbb{R}}^{n}\setminus\{0\}\,,
\end{eqnarray}
for all $l\in\{1,\dots,n\}$.  Next we note that the function $A_{2,j}(\frac{x}{|x|},0)$ is positively homogeneous of degree $0$ in the variable $x\in {\mathbb{R}}^n\setminus\{0\}$. By Proposition \ref{mldc_prop:grafun}, the function $A_{2,j}(\frac{x}{|x|},0)$ is even. Then the function $J_{jl,1}$ is odd and positively homogeneous of degree $-(n-1)$ for all $l\in\{1,\dots,n\}$ and Theorem \ref{mldc_thm:mirandao01} implies that
\begin{eqnarray}\label{mldc_prop:pdersa010om7}
&&K^+[J_{jl,1},\cdot]\in {\mathcal{L}}\left(C^{0,1}(\partial\Omega),C^{0,\omega_1(\cdot)}(\overline{\Omega})\right)\,,
\\	\nonumber
&&K^-[J_{jl,1},\cdot]\in {\mathcal{L}}\left(C^{0,1}(\partial\Omega),C^{0,\omega_1(\cdot)}(\overline{{\mathbb{B}}_n(0,r)}\setminus\Omega)\right)\,.
\end{eqnarray}
Next we note that the function $\frac{\partial A_{2,j}}{\partial\xi_s}(\frac{x}{|x|},0)$ is positively homogeneous of degree $0$ in the variable $x\in {\mathbb{R}}^n\setminus\{0\}$. By Proposition \ref{mldc_prop:grafun}, the function $\frac{\partial A_{2,j}}{\partial\xi_s}(\frac{x}{|x|},0)$ is odd. Then the function $J_{jl,2}$ is odd and positively homogeneous of degree $-(n-1)$ for all $l\in\{1,\dots,n\}$ and Theorem \ref{mldc_thm:mirandao01} implies that
\begin{eqnarray}\label{mldc_prop:pdersa010om8}
&&K^+[J_{jl,2},\cdot]\in {\mathcal{L}}\left(C^{0,1}(\partial\Omega),C^{0,\omega_1(\cdot)}(\overline{\Omega})\right)\,,
\\	\nonumber
&&K^-[J_{jl,2},\cdot]\in {\mathcal{L}}\left(C^{0,1}(\partial\Omega),C^{0,\omega_1(\cdot)}(\overline{{\mathbb{B}}_n(0,r)}\setminus\Omega)\right)\,.
\end{eqnarray}
 Since $\partial\Omega$ is upper $(n-1)$-Ahlfors regular with respect to ${\mathbb{R}}^n$ and the singularity of $J_{jl,3}(x-\cdot)$ is weak on $\partial\Omega$ for all $x\in {\mathbb{R}}^n$ for $n\geq 3$, we have
 \begin{eqnarray}\label{mldc_prop:pdersa010om9}
&&K^+[J_{jl,3},\cdot]\in {\mathcal{L}}\left(C^{0,1}(\partial\Omega),C^{0}(\overline{\Omega})\right)\,,
\\	\nonumber
&&K^-[J_{jl,3},\cdot]\in {\mathcal{L}}\left(C^{0,1}(\partial\Omega),C^{0}(\overline{{\mathbb{B}}_n(0,r)}\setminus\Omega)\right)\,,
\\	\nonumber
&& \text{for}\ n\geq 3
\,,
\end{eqnarray}
(cf.~\cite[Prop.~6.5]{La24}, \cite[Prop.~4.1]{La22a}, \cite[Prop.~4.3]{La22b}). If instead $n=2$, we note that
the function $J_{jl,3}$ is odd and positively homogeneous of degree $-1$  for all $l\in\{1,\dots,n\}$ and Theorem \ref{mldc_thm:mirandao01} implies that
\begin{eqnarray}\label{mldc_prop:pdersa010om10}
&&K^+[J_{jl,3},\cdot]\in {\mathcal{L}}\left(C^{0,1}(\partial\Omega),C^{0,\omega_1(\cdot)}(\overline{\Omega})\right)\,,
\\	\nonumber
&&K^-[J_{jl,3},\cdot]\in {\mathcal{L}}\left(C^{0,1}(\partial\Omega),C^{0,\omega_1(\cdot)}(\overline{{\mathbb{B}}_n(0,r)}\setminus\Omega)\right)
\\	\nonumber
&& \text{for}\ n=2
\,.
\end{eqnarray}
Since $\partial\Omega$ is upper $(n-1)$-Ahlfors regular with respect to ${\mathbb{R}}^n$ and the singularity of $\frac{\partial k_j}{\partial x_l}(x-\cdot)-\sum_{s=1}^3
J_{jl,s}(x-\cdot)$ is weak on $\partial\Omega$  for all $x\in {\mathbb{R}}^n$, we have
 \begin{eqnarray}\label{mldc_prop:pdersa010om11}
&&K^+\left[
\frac{\partial k_j}{\partial x_l}-\sum_{s=1}^3
J_{jl,s} 
,\cdot\right]\in {\mathcal{L}}\left(C^{0,1}(\partial\Omega),C^{0}(\overline{\Omega})\right)\,,
\\	\nonumber
&&K^-\left[\frac{\partial k_j}{\partial x_l}-\sum_{s=1}^3
J_{jl,s},\cdot\right]\in {\mathcal{L}}\left(C^{0,1}(\partial\Omega),C^{0}(\overline{{\mathbb{B}}_n(0,r)}\setminus\Omega)\right) 
\,,
\end{eqnarray}
(cf.~\cite[Prop.~6.5]{La24}, \cite[Prop.~4.1]{La22a}, \cite[Prop.~4.3]{La22b} and see  the notation in (\ref{mldc_eq:1miraslay})). By equalities (\ref{mldc_prop:pdersa010om5}), (\ref{mldc_prop:pdersa010om6}) and by the memberships of 
(\ref{mldc_prop:pdersa010om7})--(\ref{mldc_prop:pdersa010om11}), we deduce the validity of
the memberships of (\ref{mldc_prop:pdersa010om4}) and thus the proof is complete.
\hfill  $\Box$ 

\vspace{\baselineskip}

 We are now ready to prove the following extension of a known result of Miranda \cite{Mi65} on the single layer potential that is associated to the fundamental solution of elliptic operators of second order with constant coefficients. Miranda had considered the case of homogeneous operators and the case of domains of class $C^{m,\alpha}$  and of densities $\mu \in C^{m,\alpha}(\partial\Omega)$ for $\alpha\in]0,1[$. Later Dalla Riva \cite{Da13} has extended the result of Miranda by removing the restriction that the operator be homogeneous, under the assumption that $\alpha\in]0,1[$. Here, we consider the limiting case in which $\alpha=1$ and we first state and  prove case $m=1$. 
  \begin{theorem}\label{mldc_thm:sla010om}
 Let ${\mathbf{a}}$ be as in (\ref{mldc_introd0}), (\ref{mldc_ellip}), (\ref{mldc_symr}).  Let $S_{ {\mathbf{a}} }$ be a fundamental solution of $P[{\mathbf{a}},D]$.  
  Let $\Omega$ be a bounded open subset of ${\mathbb{R}}^n$ of class $C^{1,1}$.    
 Then the following statements hold.
\begin{enumerate}
\item[(i)] If  $\mu\in C^{0,1}(\partial\Omega)$, the map $v_\Omega^+\left[S_{ {\mathbf{a}} },\mu\right]$ belongs to $C^{1,\omega_1(\cdot)}(\overline{\Omega})$. Moreover, the map from $C^{0,1}(\partial\Omega)$ to $C^{1,\omega_1(\cdot)}(\overline{\Omega})$ that takes $\mu$ to $v_\Omega^+\left[S_{ {\mathbf{a}} },\mu\right]$ is linear and continuous.
\item[(ii)]  Let $r\in]0,+\infty[$ be such that $\overline{\Omega}\subseteq {\mathbb{B}}_n(0,r)$. If  $\mu\in C^{0,1}(\partial\Omega)$, the restriction  $v_\Omega^-\left[S_{ {\mathbf{a}} },\mu\right]_{|\overline{{\mathbb{B}}_n(0,r)}\setminus\Omega}$ belongs to $C^{1,\omega_1(\cdot)}(\overline{{\mathbb{B}}_n(0,r)}\setminus\Omega)$.
Moreover, the map from $C^{0,1}(\partial\Omega)$ to $C^{1,\omega_1(\cdot)}(\overline{{\mathbb{B}}_n(0,r)}\setminus\Omega)$ that takes $\mu$ to $v_\Omega^-\left[S_{ {\mathbf{a}} },\mu\right]_{|\overline{{\mathbb{B}}_n(0,r)}\setminus\Omega}$ is linear and continuous.
 
\end{enumerate}
\end{theorem}
{\bf Proof.}   Let $\gamma\in]0,1[$. By  \cite[Lemma 4.2 (i)]{DoLa17}, we have
\[
C_{0, S_{ {\mathbf{a}} },\overline{{\mathbb{B}}_n(0,2r)},n-1-\gamma}\equiv
\sup_{ 0<|\xi|\leq 4r}|\xi|^{n-1-\gamma}| S_{ {\mathbf{a}} } (\xi)|<+\infty\,.
\]
Since $\partial\Omega$ is upper $(n-1)$-Ahlfors regular with respect to ${\mathbb{R}}^n$ (cf.~\cite[(1.4)]{La22a}, \cite[Prop.~6.5]{La24}), then Lemma 3.4 of \cite{La22a} implies that
\[ 
\sup_{x\in \overline{{\mathbb{B}}_n(0,r)}}
\int_{\partial\Omega}\frac{d\sigma_y}{|x-y|^{n-1-\gamma}}
<+\infty\,,
\]
 (an inequality that one could also prove directly by elementary calculus). Then
\begin{eqnarray*}
\lefteqn{
\sup_{x\in \overline{{\mathbb{B}}_n(0,r)}}
\left|
\int_{\partial\Omega}S_{ {\mathbf{a}} }(x-y)\mu(y)\,d\sigma_y
\right|
}
\\ \nonumber
&&\qquad\qquad  
\leq C_{0, S_{ {\mathbf{a}} },\overline{{\mathbb{B}}_n(0,2r)},n-1-\gamma}
\sup_{x\in \overline{{\mathbb{B}}_n(0,r)}}
\int_{\partial\Omega}\frac{d\sigma_y}{|x-y|^{n-1-\gamma}}
\|\mu\|_{C^{0}(\partial\Omega)}
\,,
\end{eqnarray*}
for all $\mu\in C^{0,1}(\partial\Omega)$. Moreover, $v[S_{ {\mathbf{a}} },\mu]$ is continuous in ${\mathbb{R}}^n$ for all $\mu\in C^{0}(\partial\Omega)$ (cf.~\textit{e.g.}, \cite[Prop.~4.3]{La22b})) and accordingly,
\begin{eqnarray}\label{mldc_thm:sla010om1}
&&v^+[S_{ {\mathbf{a}} },\cdot]\in {\mathcal{L}}\left(C^{0,1}(\partial\Omega),C^{0}(\overline{\Omega})\right)\,,
\\	\nonumber
&&v^-[S_{ {\mathbf{a}} },\cdot]\in {\mathcal{L}}\left(C^{0,1}(\partial\Omega),C^{0}(\overline{{\mathbb{B}}_n(0,r)}\setminus\Omega)\right)\,.
\end{eqnarray}
By the cassical differentiation theorem for integrals depending on a parameter, we have 
\[
\frac{\partial}{\partial x_j}\int_{\partial\Omega}S_{ {\mathbf{a}} }(x-y)\mu(y)\,d\sigma_y
=\int_{\partial\Omega}\frac{\partial}{\partial x_j}S_{ {\mathbf{a}} }(x-y)\mu(y)\,d\sigma_y\qquad\forall x\in {\mathbb{R}}^n\setminus\partial\Omega
\]
for all $j\in\{1,\dots,n\}$,   $\mu\in C^{0}(\partial\Omega)$ and accordingly Proposition \ref{mldc_prop:pdersa010om} implies that
\begin{eqnarray}\label{mldc_thm:sla010om2}
&&\frac{\partial}{\partial x_j}v^+[S_{ {\mathbf{a}} },\cdot]\in {\mathcal{L}}\left(C^{0,1}(\partial\Omega),C^{0,\omega_1(\cdot)}(\overline{\Omega})\right)\,,
\\	\nonumber
&&\frac{\partial}{\partial x_j}v^-[S_{ {\mathbf{a}} },\cdot]\in {\mathcal{L}}\left(C^{0,1}(\partial\Omega),C^{0,\omega_1(\cdot)}(\overline{{\mathbb{B}}_n(0,r)}\setminus\Omega)\right)\,.
\end{eqnarray}
for all $j\in\{1,\dots,n\}$. Then the memberships of (\ref{mldc_thm:sla010om1}) and (\ref{mldc_thm:sla010om2}) imply the validity of statements (i), (ii) for $m=1$.

\section{A limiting case of a theorem of C.~Miranda for the double layer potential in case $m=1$}\label{mldc_sec:1miradolay}

We first introduce the following statement that follows by the definition of double layer potential and by Theorem \ref{mldc_thm:sla010om} on the single layer potential.
\begin{theorem}\label{mldc_thm:dlay010om}
 Let ${\mathbf{a}}$ be as in (\ref{mldc_introd0}), (\ref{mldc_ellip}), (\ref{mldc_symr}).  Let $S_{ {\mathbf{a}} }$ be a fundamental solution of $P[{\mathbf{a}},D]$.  
  Let $\Omega$ be a bounded open subset of ${\mathbb{R}}^n$ of class $C^{1,1}$.    
 Then the following statements hold.
\begin{enumerate}
\item[(i)] If  $\mu\in C^{0,1}(\partial\Omega)$,  then the restriction 
$w_\Omega[{\mathbf{a}},S_{{\mathbf{a}}},\mu]_{|\Omega}$ can be extended uniquely to   function $w^{+}_\Omega[ {\mathbf{a}},S_{{\mathbf{a}}},\mu]\in C^{0,\omega_1(\cdot)}(\overline{\Omega})$. Moreover, the map from $C^{0,1}(\partial\Omega)$ to $C^{0,\omega_1(\cdot)}(\overline{\Omega})$ that takes $\mu$ to $w_\Omega^+\left[{\mathbf{a}},S_{ {\mathbf{a}} },\mu\right]$ is linear and continuous.
\item[(ii)]  Let $r\in]0,+\infty[$ be such that $\overline{\Omega}\subseteq {\mathbb{B}}_n(0,r)$. If  $\mu\in C^{0,1}(\partial\Omega)$,  then the restriction 
$w_\Omega[ {\mathbf{a}},S_{{\mathbf{a}}},\mu]_{|\Omega^{-}}$ can be extended uniquely to a continuous function $w^{-}_\Omega[ {\mathbf{a}},S_{{\mathbf{a}}},\mu]$ from $\overline{\Omega^{-}}$ to ${\mathbb{C}}$ and 
the restriction  $w_\Omega^-\left[{\mathbf{a}},S_{ {\mathbf{a}} },\mu\right]_{|\overline{{\mathbb{B}}_n(0,r)}\setminus\Omega}$ belongs to $C^{m,\omega_1(\cdot)}(\overline{{\mathbb{B}}_n(0,r)}\setminus\Omega)$.
Moreover, the map from $C^{0,1}(\partial\Omega)$ to $C^{0,\omega_1(\cdot)}(\overline{{\mathbb{B}}_n(0,r)}\setminus\Omega)$ that takes $\mu$ to $w_\Omega^-\left[{\mathbf{a}},S_{ {\mathbf{a}} },\mu\right]_{|\overline{{\mathbb{B}}_n(0,r)}\setminus\Omega}$ is linear and continuous.
 \end{enumerate}
\end{theorem}
{\bf Proof.} (i) Formula (\ref{mldc_introd3}) implies that 
\begin{eqnarray}\label{mldc_thm:dlay010om1}
\lefteqn{
w_\Omega^+[{\mathbf{a}},S_{{\mathbf{a}}}   ,\mu](x)
}
\\ \nonumber
&&\qquad
=-\sum_{l,j=1}^{n} a_{jl}\frac{\partial  }{\partial x_{j}}
v_\Omega^+[S_{{\mathbf{a}}},\mu \nu_{l}](x)
-
\sum_{l=1}^{n}  a_{l}
v_\Omega^+[S_{{\mathbf{a}}},\mu\nu_{l}](x)\qquad\forall x\in\Omega\,.
 \end{eqnarray}
 Since the components of $\nu$ are of class $C^{0,1}(\partial\Omega)$ and the pointwise product is continuous in $C^{0,1}(\partial\Omega)$, then
 Then Proposition \ref{mldc_prop:pdersa010om},   Theorem   \ref{mldc_thm:sla010om} 
 and the continuity of the embedding of $C^{1,\omega_1(\cdot)}(\overline{\Omega})$ into $C^{0,\omega_1(\cdot)}(\overline{\Omega})$ 
  imply that $w_\Omega^+[{\mathbf{a}},S_{{\mathbf{a}}}   ,\mu]_{|\Omega}$ can be extended to a unique element of $C^{0,\omega_1(\cdot)}(\overline{\Omega})$ and that statement (i) holds true. Then statement (ii) can be deduced applying statement (i) to the open bounded set $\mathbb{B}_n(0,r)\setminus\overline{\Omega}$. Indeed the map from $C^{0,1} (\partial\Omega)$ to $C^{0,1}(\partial(\mathbb{B}_n(0,r)\setminus\overline{\Omega}))$ that  takes $\mu$ to
\begin{equation}\label{mldc_thm:dlay010om2}
\tilde\mu(x)\equiv
\begin{cases}
\mu(x)&\text{if }x\in\partial\Omega\,,\\
0&\text{if }x\in\partial\mathbb{B}_n(0,r)
\end{cases}
\end{equation}
is linear and continuous and we have
\begin{equation}\label{mldc_thm:dlay010om3}
w_\Omega^-[{\mathbf{a}},S_{{\mathbf{a}}}   ,\mu]_{|\overline{\mathbb{B}_n(0,r)}\setminus\Omega}
=-w_{	
\mathbb{B}_n(0,r)\setminus\overline{\Omega}	}^+[{\mathbf{a}},S_{{\mathbf{a}}}   ,\tilde\mu]
\end{equation}
for all $\mu  \in C^{0,1}(\partial\Omega)$.\hfill  $\Box$ 

\vspace{\baselineskip}

In order proceed further, we need the definition of tangential derivative   that we now introduce.
Let $\Omega$ be an open subset of ${\mathbb{R}}^n$ of class $C^1$.  If $l,r\in\{1,\dots,n\}$,  then $M_{lr}$ denotes the tangential derivative 
 operator from $C^{1}(\partial\Omega)$ to $C^{0}(\partial\Omega)$ that takes $f$ to  
 \[
M_{lr}[f]\equiv \nu_{l}\frac{\partial f^\sharp}{\partial x_{r}}-
\nu_{r}\frac{\partial f^\sharp}{\partial x_{l}}\qquad {\text{on}}\ \partial\Omega\,,
\]
where  $f^\sharp$ is any continuously differentiable  extension of $f$  to an open neighborhood of $\partial\Omega$. We note that $M_{lr}[f]$ is independent of the specific choice of $f^\sharp$ (cf.~\textit{e.g.}, Dalla Riva, the author and Musolino \cite[\S 2.21]{DaLaMu21}). If necessary, we write
$M_{lr,x}$ to emphasize that we are taking $x$ as variable of the differential operator $M_{lr}$. 
Then by exploiting standard differentiation rules, one can prove the following statement (cf.~\textit{e.g.} Dondi and the author \cite[(7.1)]{DoLa17})
  \begin{proposition}\label{mldc_prop:dlayder}
 Let ${\mathbf{a}}$ be as in (\ref{mldc_introd0}), (\ref{mldc_ellip}), (\ref{mldc_symr}).  Let $S_{ {\mathbf{a}} }$ be a fundamental solution of $P[{\mathbf{a}},D]$.  Let $\alpha\in]0,1]$.
   Let $\Omega$ be a bounded open subset of ${\mathbb{R}}^n$ of class $C^{1,\alpha}$.    Let $j\in\{1,\dots,n\}$. If $\mu\in C^{1,\alpha}(\partial\Omega)$, then the following equality holds
\begin{eqnarray}
\label{mldc_prop:dlayder0}
\lefteqn{\frac{\partial}{\partial x_{j}}w_\Omega[{\mathbf{a}},S_{{\mathbf{a}}},\mu](x)}
\\   \nonumber
&&
=\sum_{s,l=1}^{n}a_{ls} \frac{\partial}{\partial x_{l}}
\biggl\{
\int_{\partial\Omega}
S_{{\mathbf{a}}}(x-y)
M_{js}[\mu](y)\,d\sigma_{y}
\biggr\}
\\   \nonumber
&&\ \ +
\int_{\partial\Omega}
\biggl[
 DS_{{\mathbf{a}}}(x-y)a^{(1)}+
a S_{{\mathbf{a}}}(x-y)
\biggr]   \nu_{j} (y) \mu(y)\,d\sigma_{y}
\\   \nonumber
&\ &\ \ -\int_{\partial\Omega}\partial_{x_{j}}S_{{\mathbf{a}}}(x-y) 
\nu^{t}(y)a^{(1)}\mu(y)\,d\sigma_{y} \qquad\forall x\in{\mathbb{R}}^{n}\setminus\partial\Omega\,.
\end{eqnarray}
\end{proposition}
We note that formula (\ref{mldc_prop:dlayder0}) for the Laplace operator with $n=3$ can be found in 
G\"{u}nter~\cite[Ch.~2, \S\ 10, (42)]{Gu67}. By combining Theorems \ref{mldc_thm:sla010om} and Proposition \ref{mldc_prop:dlayder}, we deduce that under the assumptions of Proposition \ref{mldc_prop:dlayder}, the following equality holds
\begin{eqnarray}
\label{mldc_prop:dlayder1}
\lefteqn{\frac{\partial}{\partial x_{j}}w^{+}_\Omega[{\mathbf{a}},S_{{\mathbf{a}}},\mu] }
\\  \nonumber
&&
=\sum_{s,l=1}^{n}a_{ls} \frac{\partial}{\partial x_{l}}
v^{+}_\Omega[S_{{\mathbf{a}}},M_{js}[\mu]] 
+
Dv^{+}_\Omega[S_{{\mathbf{a}}},\nu_{j}   \mu] a^{(1)}
\\    \nonumber
&&\ \ 
+av^{+}_\Omega[S_{{\mathbf{a}}},\nu_{j}   \mu] 
 -\frac{\partial}{\partial x_{j}}v^{+}_\Omega[S_{{\mathbf{a}}},
 (\nu^{t} a^{(1)})\mu]
\qquad{\mathrm{on}}\ \overline{\Omega}\,.
\end{eqnarray}
We are now ready to prove the following extension of a Theorem of Miranda for the double layer potential to the limiting case in which $\alpha=1$.
 \begin{theorem}\label{mldc_thm:dlay111om}
 Let ${\mathbf{a}}$ be as in (\ref{mldc_introd0}), (\ref{mldc_ellip}), (\ref{mldc_symr}).  Let $S_{ {\mathbf{a}} }$ be a fundamental solution of $P[{\mathbf{a}},D]$.  
  Let $\Omega$ be a bounded open subset of ${\mathbb{R}}^n$ of class $C^{1,1}$.    
 Then the following statements hold.
\begin{enumerate}
\item[(i)] If  $\mu\in C^{1,1}(\partial\Omega)$, the map $w_\Omega^+\left[{\mathbf{a}},S_{ {\mathbf{a}} },\mu\right]$ belongs to $C^{1,\omega_1(\cdot)}(\overline{\Omega})$. Moreover, the map from $C^{1,1}(\partial\Omega)$ to $C^{1,\omega_1(\cdot)}(\overline{\Omega})$ that takes $\mu$ to $w_\Omega^+\left[{\mathbf{a}},S_{ {\mathbf{a}} },\mu\right]$ is linear and continuous.
\item[(ii)]  Let $r\in]0,+\infty[$ be such that $\overline{\Omega}\subseteq {\mathbb{B}}_n(0,r)$. If  $\mu\in C^{1,1}(\partial\Omega)$, the restriction  $w_\Omega^-\left[{\mathbf{a}},S_{ {\mathbf{a}} },\mu\right]_{|\overline{{\mathbb{B}}_n(0,r)}\setminus\Omega}$ belongs to $C^{1,\omega_1(\cdot)}(\overline{{\mathbb{B}}_n(0,r)}\setminus\Omega)$.
Moreover, the map from $C^{1,1}(\partial\Omega)$ to $C^{1,\omega_1(\cdot)}(\overline{{\mathbb{B}}_n(0,r)}\setminus\Omega)$ that takes $\mu$ to $w_\Omega^-\left[{\mathbf{a}},S_{ {\mathbf{a}} },\mu\right]_{|\overline{{\mathbb{B}}_n(0,r)}\setminus\Omega}$ is linear and continuous.
 \end{enumerate}
\end{theorem}
 {\bf Proof.} (i) By the continuity of the embedding of $C^{1,1}(\partial\Omega)$ into $C^{0,1}(\partial\Omega)$ and by Theorem \ref{mldc_thm:dlay010om}, we have
 \[
w_\Omega^+\left[{\mathbf{a}},S_{ {\mathbf{a}} },\cdot\right]
\in 
{\mathcal{L}}\left(C^{1,1}(\partial\Omega),C^{0,\omega_1(\cdot)}(\overline{\Omega})\right)
\]
Thus it suffices to show that
\begin{equation}\label{mldc_thm:dlay111om2}
\frac{\partial}{\partial x_j}w_\Omega^+\left[{\mathbf{a}},S_{ {\mathbf{a}} },\cdot\right]
\in 
{\mathcal{L}}\left(C^{1,1}(\partial\Omega),C^{0,\omega_1(\cdot)}(\overline{\Omega})\right)
\end{equation}
for all $j\in\{1,\dots,n\}$.  We plan to do so by exploiting formula  (\ref{mldc_prop:dlayder1}) and Theorem \ref{mldc_thm:sla010om} on the single layer potential. Since the tangential derivative $M_{js}[\cdot]$ is linear and continuous from  $C^{1,1}(\partial\Omega)$ to $C^{0,1}(\partial\Omega)$, Theorem   \ref{mldc_thm:sla010om} implies that the first addendum in the right hand side of  formula  (\ref{mldc_prop:dlayder1})  is linear and continuous from $C^{1,1}(\partial\Omega)$ to $C^{0,\omega_1(\cdot)}(\overline{\Omega})$. Since the components of $\nu$ are of class $C^{0,1}(\partial\Omega)$, $C^{1,1}(\partial\Omega)$ is continuosly embedded into $C^{0,1}(\partial\Omega)$ and the pointwise product is continuous in $C^{0,1}(\partial\Omega)$
 Theorem   \ref{mldc_thm:sla010om} implies that the second and fourth addendum in the right hand side of  formula  (\ref{mldc_prop:dlayder1})  define linear and continuous maps from $C^{1,1}(\partial\Omega)$ to $C^{0,\omega_1(\cdot)}(\overline{\Omega})$. Since the components of $\nu$ are of class $C^{0,1}(\partial\Omega)$,  $C^{1,1}(\partial\Omega)$ is continuosly embedded into $C^{0,1}(\partial\Omega)$  and the pointwise product is continuous in $C^{0,1}(\partial\Omega)$
 Theorem   \ref{mldc_thm:sla010om} and the continuity of the embedding of $C^{1,\omega_1(\cdot)}(\overline{\Omega})$ into $C^{0,\omega_1(\cdot)}(\overline{\Omega})$ imply that the third addendum in the right hand side of  formula  (\ref{mldc_prop:dlayder1})  defines a linear and continuous map from $C^{1,1}(\partial\Omega)$ to $C^{0,\omega_1(\cdot)}(\overline{\Omega})$. Hence, the membership of (\ref{mldc_thm:dlay111om2}) holds  true and the proof of (i) is complete.

 Then statement (ii) can be deduced applying statement (i) to the open bounded set $\mathbb{B}_n(0,r)\setminus\overline{\Omega}$ by exploiting  (\ref{mldc_thm:dlay010om2}) and (\ref{mldc_thm:dlay010om3}).\hfill  $\Box$ 

\vspace{\baselineskip}

\section{A limiting case of a theorem of C.~Miranda for the single and double layer potential in case $m\geq1$}\label{mldc_sec:miradslay}

By exploiting standard differentiation rules, one can prove the following statement. 
  \begin{proposition}\label{mldc_prop:slayder}
 Let ${\mathbf{a}}$ be as in (\ref{mldc_introd0}), (\ref{mldc_ellip}), (\ref{mldc_symr}).  Let $S_{ {\mathbf{a}} }$ be a fundamental solution of $P[{\mathbf{a}},D]$.   
   Let $\Omega$ be a bounded open subset of ${\mathbb{R}}^n$ of class $C^{2}$.    Let $j\in\{1,\dots,n\}$. Then the following equality holds
\begin{eqnarray}
\label{mldc_prop:slayder0}
\lefteqn{\frac{\partial}{\partial x_{j}}v_\Omega[S_{{\mathbf{a}}},\mu](x)
=v_\Omega
\left[S_{{\mathbf{a}}},
\sum_{r,s=1}^{n}M_{rj}\left[\frac{a_{rs}\nu_s\mu}{\nu^ta^{(2)}\nu}\right]\right](x)
}
\\ \nonumber
&&\qquad 
-w_\Omega\left[{\mathbf{a}},S_{{\mathbf{a}}},\frac{\nu_j\mu}{\nu^ta^{(2)}\nu}\right](x)
-v_\Omega\left[S_{{\mathbf{a}}},\frac{(a^{(1)}\nu)\nu_j\mu}{\nu^ta^{(2)}\nu}\right](x)
\quad\forall x\in{\mathbb{R}}^{n}\setminus\partial\Omega\,,
\end{eqnarray}
for all $\mu\in C^{1}(\partial\Omega)$. 
\end{proposition}
{\bf Proof.} If $x\in{\mathbb{R}}^{n}\setminus\partial\Omega$, then the classical differentiation theorem for integrals depending on  a parameter   implies that
\begin{eqnarray*}
\lefteqn{
\frac{\partial}{\partial x_{j}}v_\Omega[S_{{\mathbf{a}}},\mu](x)
}
\\ \nonumber
&&\qquad\qquad
=\sum_{r,s=1}^n\int_{\partial\Omega}
-\frac{\partial}{\partial y_j}
\left(S_{{\mathbf{a}}}(x-y)\right)
\nu_r(y)a_{rs}\nu_s(y)\mu(y)\frac{d\sigma_y}{\nu^t(y)a^{(2)}\nu(y)}
\\ \nonumber
&&\qquad\qquad\quad
+\sum_{r,s=1}^n\int_{\partial\Omega}\frac{\partial}{\partial y_r}
\left(S_{{\mathbf{a}}}(x-y)\right)\nu_j(y)a_{rs}\nu_s(y)\mu(y)\frac{d\sigma_y}{\nu^t(y)a^{(2)}\nu(y)}
\\ \nonumber
&&\qquad\qquad\quad
-\sum_{r,s=1}^n\int_{\partial\Omega}\frac{\partial}{\partial y_r}
\left(S_{{\mathbf{a}}}(x-y)\right)\nu_j(y)a_{rs}\nu_s(y)\mu(y)\frac{d\sigma_y}{\nu^t(y)a^{(2)}\nu(y)}
\\ \nonumber
&&\qquad\qquad
=-\sum_{r,s=1}^n\int_{\partial\Omega}M_{rj,y}[S_{{\mathbf{a}}}(x-y)]a_{rs}\nu_s(y)\mu(y)\frac{d\sigma_y}{\nu^t(y)a^{(2)}\nu(y)}
\\ \nonumber
&&\qquad\qquad\quad
+\sum_{r,s=1}^n\int_{\partial\Omega}\nu_j(y)\mu(y)a_{rs}\nu_s(y)\frac{\partial S_{{\mathbf{a}}}}{\partial x_r}(x-y)\frac{d\sigma_y}{\nu^t(y)a^{(2)}\nu(y)}
\\ \nonumber
&&\qquad\qquad\quad
+\int_{\partial\Omega}\nu_j(y)\mu(y)\sum_{l=1}^n\nu_l(y)a_lS_{{\mathbf{a}}}(x-y)\frac{d\sigma_y}{\nu^t(y)a^{(2)}\nu(y)}
\\ \nonumber
&&\qquad\qquad\quad
-\int_{\partial\Omega}\nu_j(y)\mu(y)\sum_{l=1}^n\nu_l(y)a_lS_{{\mathbf{a}}}(x-y)\frac{d\sigma_y}{\nu^t(y)a^{(2)}\nu(y)}
\\ \nonumber
&&\qquad\qquad
=\sum_{r,s=1}^n\int_{\partial\Omega} S_{{\mathbf{a}}}(x-y)M_{rj}\left[\frac{a_{rs}\nu_s \mu}{ \nu^ta^{(2)}\nu}\right](y) d\sigma_y 
\\ \nonumber
&&\qquad\qquad\quad
-w_\Omega\left[{\mathbf{a}},S_{{\mathbf{a}}},\frac{\nu_j\mu}{\nu^ta^{(2)}\nu}\right](x)
-v_\Omega\left[S_{{\mathbf{a}}},\frac{(a^{(1)}\nu)\nu_j\mu}{\nu^ta^{(2)}\nu}\right](x)
\end{eqnarray*}
and thus equality (\ref{mldc_prop:slayder0}) holds true (see also   Lemma 2.86 of \cite{DaLaMu21}  on the tangential derivative).\hfill  $\Box$ 

\vspace{\baselineskip}

We note that a related form of formula (\ref{mldc_prop:slayder0}) for the Laplace operator with $n=3$ can be found in 
G\"{u}nter~\cite[Ch.~2, \S\ 8, (38)]{Gu67}. We are now ready to prove the following extension of a known result of Miranda \cite{Mi65} on the single and double layer potential associated to the fundamental solution of elliptic operators of second order with constant coefficients. Miranda had considered the case of homogeneous operators and the case of domains of class $C^{m,\alpha}$  and of densities $\mu \in C^{m-1,\alpha}(\partial\Omega)$ for the single layer potential
and $\mu \in C^{m,\alpha}(\partial\Omega)$ for the double layer potential when $\alpha\in]0,1[$. Later Dalla Riva \cite{Da13} has extended the result of Miranda by removing the restriction that the operator be homogeneous, under the assumption that $\alpha\in]0,1[$. Here, we consider the limiting case in which   $\alpha=1$. 
\begin{theorem}\label{mldc_thm:dslam1mom}
Let ${\mathbf{a}}$ be as in (\ref{mldc_introd0}), (\ref{mldc_ellip}), (\ref{mldc_symr}).  Let $S_{ {\mathbf{a}} }$ be a fundamental solution of $P[{\mathbf{a}},D]$.  
 Let $m\in {\mathbb{N}}\setminus\{0\}$. Let $\Omega$ be a bounded open subset of ${\mathbb{R}}^n$ of class $C^{m,1}$.    
 Then the following statements hold.
\begin{enumerate}
\item[(i)] If  $\mu\in C^{m-1,1}(\partial\Omega)$, the map $v_\Omega^+\left[S_{ {\mathbf{a}} },\mu\right]$ belongs to $C^{m,\omega_1(\cdot)}(\overline{\Omega})$. Moreover, the map from $C^{m-1,1}(\partial\Omega)$ to $C^{m,\omega_1(\cdot)}(\overline{\Omega})$ that takes $\mu$ to $v_\Omega^+\left[S_{ {\mathbf{a}} },\mu\right]$ is linear and continuous.
\item[(ii)]  Let $r\in]0,+\infty[$ be such that $\overline{\Omega}\subseteq {\mathbb{B}}_n(0,r)$. If  $\mu\in C^{m-1,1}(\partial\Omega)$, the restriction  $v_\Omega^-\left[S_{ {\mathbf{a}} },\mu\right]_{|\overline{{\mathbb{B}}_n(0,r)}\setminus\Omega}$ belongs to $C^{m,\omega_1(\cdot)}(\overline{{\mathbb{B}}_n(0,r)}\setminus\Omega)$.
Moreover, the map from $C^{m-1,1}(\partial\Omega)$ to $C^{m,\omega_1(\cdot)}(\overline{{\mathbb{B}}_n(0,r)}\setminus\Omega)$ that takes $\mu$ to $v_\Omega^-\left[S_{ {\mathbf{a}} },\mu\right]_{|\overline{{\mathbb{B}}_n(0,r)}\setminus\Omega}$ is linear and continuous.
\item[(iii)] If  $\mu\in C^{m,1}(\partial\Omega)$, the map $w_\Omega^+\left[{\mathbf{a}},S_{ {\mathbf{a}} },\mu\right]$ belongs to $C^{m,\omega_1(\cdot)}(\overline{\Omega})$. Moreover, the map from $C^{m,1}(\partial\Omega)$ to $C^{m,\omega_1(\cdot)}(\overline{\Omega})$ that takes $\mu$ to $w_\Omega^+\left[{\mathbf{a}},S_{ {\mathbf{a}} },\mu\right]$ is linear and continuous.
\item[(iv)]  Let $r\in]0,+\infty[$ be such that $\overline{\Omega}\subseteq {\mathbb{B}}_n(0,r)$. If  $\mu\in C^{m,1}(\partial\Omega)$, the restriction  $w_\Omega^-\left[{\mathbf{a}},S_{ {\mathbf{a}} },\mu\right]_{|\overline{{\mathbb{B}}_n(0,r)}\setminus\Omega}$ belongs to $C^{m,\omega_1(\cdot)}(\overline{{\mathbb{B}}_n(0,r)}\setminus\Omega)$.
Moreover, the map from $C^{m,1}(\partial\Omega)$ to $C^{m,\omega_1(\cdot)}(\overline{{\mathbb{B}}_n(0,r)}\setminus\Omega)$ that takes $\mu$ to $w_\Omega^-\left[{\mathbf{a}},S_{ {\mathbf{a}} },\mu\right]_{|\overline{{\mathbb{B}}_n(0,r)}\setminus\Omega}$ is linear and continuous.
\end{enumerate}
\end{theorem}
{\bf Proof.}   We first prove (i), (iii)   by induction on $m$. If $m=1$, then Theorems \ref{mldc_thm:sla010om}, \ref{mldc_thm:dlay111om} imply the validity of statements (i), (iii). 
We now assume that statements (i), (iii) hold for $m\geq 1$ and we prove them for $m+1$. So we now assume that $\Omega$ is of class $C^{m+1,1}$. By case $m=1$ and by the continuity of the embedding of $C^{(m+1)-1,1}(\partial\Omega)$ into $C^{0,1}(\partial\Omega)$ and of 
  $C^{m+1,1}(\partial\Omega)$ into $C^{1,1}(\partial\Omega)$ and of $C^{1,\omega_1(\cdot)}(\overline{\Omega})$ into $C^{0}(\overline{\Omega})$, we have
\begin{eqnarray*}
&&v_\Omega^+\left[S_{ {\mathbf{a}} },\cdot\right]
\in {\mathcal{L}}\left(C^{(m+1)-1,1}(\partial\Omega),C^{0}(\overline{\Omega})\right)\,,
\\ \nonumber
&&w_\Omega^+\left[{\mathbf{a}},S_{ {\mathbf{a}} },\cdot\right]
\in {\mathcal{L}}\left(C^{m+1,1}(\partial\Omega),C^{0}(\overline{\Omega})\right)\,.
\end{eqnarray*}
Thus it suffices to show that
\begin{eqnarray}\label{mldc_thm:dslam1mom2}
&&\frac{\partial}{\partial x_j}v_\Omega^+\left[S_{ {\mathbf{a}} },\cdot\right]
\in {\mathcal{L}}\left(C^{(m+1)-1,1}(\partial\Omega),C^{m,\omega_1(\cdot)}(\overline{\Omega})\right)\,,
\\ \label{mldc_thm:dslam1mom3}
&&\frac{\partial}{\partial x_j}w_\Omega^+\left[{\mathbf{a}},S_{ {\mathbf{a}} },\cdot\right]
\in {\mathcal{L}}\left(C^{m+1,1}(\partial\Omega),C^{m,\omega_1(\cdot)}(\overline{\Omega})\right)\,,
\end{eqnarray}
for all $j\in\{1,\dots,n\}$. We plan to do so by exploiting formulas  (\ref{mldc_prop:dlayder1}) 
and (\ref{mldc_prop:slayder0}). We first consider the membership of (\ref{mldc_thm:dslam1mom2}). Since the components of $\nu$ are of class $C^{m,1}(\partial\Omega)$, the pointwise product is continuous in $C^{(m+1)-1,1}(\partial\Omega)$,   the function $(\nu^ta^{(2)}\nu)^{-1}$ belongs to $C^{m,1}(\partial\Omega)$, the tangential derivative $M_{rj}[\cdot]$ is linear and continuous from  $C^{(m+1)-1,1}(\partial\Omega)$ to $C^{m-1,1}(\partial\Omega)$, the inductive assumption of statement (i) implies that the first addendum in the right hand side of  formula (\ref{mldc_prop:slayder0}) defines a linear and continuous operator from  $C^{(m+1)-1,1}(\partial\Omega)$ to
$C^{m,\omega_1(\cdot)}(\overline{\Omega})$.

Since the components of $\nu$ are of class $C^{m,1}(\partial\Omega)$, 
the function $(\nu^ta^{(2)}\nu)^{-1}$ belongs to $C^{m,1}(\partial\Omega)$,
the pointwise product is continuous in $C^{(m+1)-1,1}(\partial\Omega)$,  
 the inductive assumption of statement (iii) implies that the second addendum in the right hand side of  formula (\ref{mldc_prop:slayder0}) defines a linear and continuous operator from  $C^{(m+1)-1,1}(\partial\Omega)$ to
$C^{m,\omega_1(\cdot)}(\overline{\Omega})$.

Since the components of $\nu$ are of class $C^{m,1}(\partial\Omega)$, 
the function $(\nu^ta^{(2)}\nu)^{-1}$ belongs to $C^{m,1}(\partial\Omega)$, 
the pointwise product is continuous in $C^{(m+1)-1,1}(\partial\Omega)$, 
 $C^{(m+1)-1,1}(\partial\Omega)$ is continuously embedded into $C^{m-1,1}(\partial\Omega)$, the inductive assumption of statement (i) implies that the third addendum in the right hand side of  formula (\ref{mldc_prop:slayder0}) defines a linear and continuous operator from  $C^{(m+1)-1,1}(\partial\Omega)$ to
$C^{m,\omega_1(\cdot)}(\overline{\Omega})$. Hence, the membership of (\ref{mldc_thm:dslam1mom2}) holds true.

We now consider the membership of (\ref{mldc_thm:dslam1mom3}). Since  the tangential derivative $M_{js}[\cdot]$ is linear and continuous from  $C^{m+1,1}(\partial\Omega)$ to $C^{m,1}(\partial\Omega)$, statement (i) for case $(m+1)$ that we have just proved above to be a consequence of the inductive assumption implies that the first addendum in the right hand side of  formula (\ref{mldc_prop:dlayder1}) defines a linear and continuous operator from  $C^{m+1,1}(\partial\Omega)$ to
$C^{m,\omega_1(\cdot)}(\overline{\Omega})$.

Since the components of $\nu$ are of class $C^{(m+1)-1,1}(\partial\Omega)$,  $C^{m+1,1}(\partial\Omega)$ is continuously embedded into $C^{m,1}(\partial\Omega)$, 
 the pointwise product is continuous in $C^{m,1}(\partial\Omega)$,   statement (i) for case $(m+1)$ that we have just proved above to be a consequence of the inductive assumption implies that the second addendum in the right hand side of  formula (\ref{mldc_prop:dlayder1}) defines a linear and continuous operator from  $C^{m+1,1}(\partial\Omega)$ to
$C^{m,\omega_1(\cdot)}(\overline{\Omega})$.

Since the components of $\nu$ are of class $C^{(m+1)-1,1}(\partial\Omega)$,  $C^{m+1,1}(\partial\Omega)$ is continuously embedded into $C^{m,1}(\partial\Omega)$, 
 the pointwise product is continuous in $C^{m,1}(\partial\Omega)$ and $C^{m,1}(\partial\Omega)$ is continuously embedded into $C^{m-1,1}(\partial\Omega)$, the inductive assumption of statement (i) implies that the third  addendum in the right hand side of  formula (\ref{mldc_prop:dlayder1}) defines a linear and continuous operator from  $C^{m+1,1}(\partial\Omega)$ to
$C^{m,\omega_1(\cdot)}(\overline{\Omega})$.

Since the components of $\nu$ are of class $C^{(m+1)-1,1}(\partial\Omega)$,  $C^{m+1,1}(\partial\Omega)$ is continuously embedded into $C^{m,1}(\partial\Omega)$
and the pointwise product is continuous in $C^{m,1}(\partial\Omega)$, statement (i) for case $(m+1)$ that we have just proved above to be a consequence of the inductive assumption implies that the fourth  addendum in the right hand side of  formula (\ref{mldc_prop:dlayder1}) defines a linear and continuous operator from  $C^{m+1,1}(\partial\Omega)$ to
$C^{m,\omega_1(\cdot)}(\overline{\Omega})$. Hence, the membership of (\ref{mldc_thm:dslam1mom3}) holds true and thus the proof of statements (i), (iii)  is complete by the induction principle.

Then statements (ii), (iv) can be deduced applying statements (i), (iii) to the open bounded set $\mathbb{B}_n(0,r)\setminus\overline{\Omega}$ by exploiting  (\ref{mldc_thm:dlay010om2}) and (\ref{mldc_thm:dlay010om3}).\hfill  $\Box$ 

\vspace{\baselineskip}

 \noindent
{\bf Statements and Declarations}\\

 \noindent
{\bf Competing interests:} This paper does not have any  conflict of interest or competing interest.

 \noindent
{\bf Acknowledgement.} 


The author  acknowledges  the support of the Research  Project GNAMPA-INdAM   $\text{CUP}\_$E53C22001930001 `Operatori differenziali e integrali in geometria spettrale' and   of the Project funded by the European Union – Next Generation EU under the National Recovery and Resilience Plan (NRRP), Mission 4 Component 2 Investment 1.1 - Call for tender PRIN 2022 No. 104 of February, 2 2022 of Italian Ministry of University and Research; Project code: 2022SENJZ3 (subject area: PE - Physical Sciences and Engineering) ``Perturbation problems and asymptotics for elliptic differential equations: variational and potential theoretic methods''.

\end{document}